\documentclass[draft]{amsart}

\newcommand{\sgn}{\operatorname{sgn}}
\newcommand{\rot}{\operatorname{rot}}

\newcommand{\R}{\boldsymbol{R}}

\newcommand{\inner}[2]{\left\langle{#1},{#2}\right\rangle}

\newcommand{\pt}[1]{\mathsf{#1}}
\renewcommand{\phi}{\varphi}

\newtheorem{theorem}{Theorem}[section]
\newtheorem{lemma}[theorem]{Lemma}
\newtheorem{e-proposition}[theorem]{Proposition}
\newtheorem{corollary}[theorem]{Corollary}
\newtheorem{e-definition}[theorem]{Definition}

\begin{document}
\title[]{%
Singularities of Blaschke normal maps\\
 of  convex surfaces
}
\author[saji]{Kentaro Saji}
\email{ksaji@gifu-u.ac.jp}
\author[umehara]{Masaaki Umehara}
\email{umehara@math.sci.osaka-u.ac.jp}
\author[yamada]{Kotaro Yamada}
\email{kotaro@math.titech.ac.jp}

\address[saji]{%
 Department of Mathematics,
 Faculty of Education,
 Gifu University,  Yanagido 1-1, Gifu 501-1193, Japan
}
\address[umehara]{%
   Department of Mathematics, Graduate School of Science,
   Osaka University,
   Toyonaka, Osaka 560-0043,
   Japan
}
\address[yamada]{%
   Department of Mathematics,
   Tokyo Institute of Technology,
   O-okayama, Meguro, Tokyo 152-8551, Japan
}
%

\maketitle

\begin{abstract}
We prove that the difference between
the numbers of positive swallowtails and negative swallowtails 
of the Blaschke normal map for a given convex surface
in affine space is equal to the Euler number 
of the subset where the affine shape operator has 
negative determinant.
\end{abstract}
\section{Introduction.}
Throughout this Note, we assume that
$M^2$ is a compact oriented $2$-manifold without boundary.
Let $\phi$ be a bundle homomorphism
of the tangent bundle $TM^2$ into a vector bundle $E$ of rank $2$ over
$M^2$. 
A point $p$ on $M^2$ is called a {\it singular point\/} if the linear
map $\phi_p:T_pM\to E_p$ is not bijective.
We denote by $\Sigma _\phi$ the set of singular points of $\phi$.
We assume that $E$ is orientable, that is,
there is a non-vanishing section
$\mu:M^2\to E^*\wedge E^*$,
where $E^*$ is the dual vector bundle of $E$.
We now fix  a metric $\inner{~}{~}$  on $E$.
Multiplying a suitable $C^\infty$-function on $M^2$, 
we may assume that $\mu(e_1,e_2)=1$ holds for any
oriented orthonormal frame $e_1$, $e_2$ on $E$.
By using a positively oriented local coordinate system $(U;u,v)$,
the {\em signed area form\/} $d\hat A$,
the signed area density function $\lambda$, and the 
({\em un-signed}) {\em area form\/} $dA$  are defined by
\begin{align*}
    d\hat A := \phi^*\mu = \lambda\,du\wedge dv,\qquad
     dA     := |\lambda|\,du\wedge dv.
\end{align*}
Both $d\hat A$  and $dA$ are independent of the choice of  $(u,v)$, and
are $2$-forms globally defined  on $M^2$.
When $\phi$ has no singular points, 
these two forms coincide up to sign.
We set
\[
  M^+:=\bigl\{p\in M^2\setminus \Sigma_\phi \,;\, d\hat A_p=dA_p\bigr\},
  \qquad
  M^-:=\bigl\{p\in M^2\setminus \Sigma_\phi\,;\, d\hat A_p=-dA_p
  \bigr\}.
\]
The singular set $\Sigma_\phi$ coincides with $\partial M^+=\partial M^-$.
A singular point $p(\in \Sigma_\phi)$ on $M^2$ is called
{\it non-degenerate\/} if the derivative $d\lambda$ does not vanish at
$p$.
In a neighborhood of a non-degenerate singular point,
the singular set can be parametrized as a regular curve $\gamma(t)$ on
$M^2$, 
called the {\em singular curve}.
The tangential direction of $\gamma$ is called the 
{\em singular direction}.
The direction of the kernel of $\phi_p$ is called the 
{\em  null direction}, which is one dimensional.
There exists a smooth non-vanishing 
vector field $\eta(t)$  along $\gamma$ pointing in the null direction,
called the {\it null vector field}.
\begin{e-definition}\label{def:ak-point}\rm
 Take a non-degenerate singular point $p\in M^2$ and 
 let $\gamma(t)$ be the singular curve satisfying  $\gamma(0)=p$.
 Then $p$ is called an {\em $A_2$-point\/}  if
 the null direction $\eta(0)$ is transversal to
 the singular direction $\dot\gamma(0)=d\gamma/dt|_{t=0}$.
 If $p$ is not an $A_2$-point, but satisfies that
 $d(\dot\gamma(t)\wedge \eta(t))/dt$
 does not vanish at $p=\gamma(0)$,
 it is called an {\em $A_3$-point\/}, where $\wedge$ is the
 exterior product on $TM^2$.
 We fix an $A_3$-point $p$.
 If the angle of the region $M^-$ {\rm(}resp. $M^+${\rm)}
 at $p$ with respect to the
 pull-back metric $ds^2:=\phi^*\inner{~}{~}$ is
 zero, then it is called a {\it positive} (resp. {\it negative}) 
 $A_3$-point.
 ($A_3$-points are either positive or negative, see \cite{SUY2}). 
\end{e-definition}
\medskip
We now fix a metric connection $D$ of $(E,\inner{~}{~})$.
Let $\gamma(t)$ be a regular curve on $M^2$ consisting only of
$A_2$-points.
Take a null vector field $\eta(t)$ such that $(\dot \gamma,\eta)$ is a
positive frame of $TM^2$ along $\gamma$.
Then 
\begin{equation}\label{eq:sing-curve-general}
 \kappa_s(t) = \sgn\bigl(d\lambda(\eta(t))\bigr)
          \frac{
           \mu\bigl(\phi(\dot\gamma(t)),D_t\phi(\dot\gamma(t))\bigr)}{
       \inner{\phi(\dot\gamma(t))}{\phi(\dot\gamma(t))}^{3/2}}
\end{equation}
is called the {\it singular curvature\/} of $\gamma$ at $t$
(see \cite{SUY1} and \cite{SUY2}).

For an oriented orthonormal frame field $e_1$ ,$e_2$
of $E$ defined on $U\subset M^2$,
there is a unique $1$-form
$\omega$ on $U$ such that
$D_Xe_1=-\omega(X)e_2$, $D_Xe_2=\omega(X)e_1$.
Then $d\omega$ does not depend on the choice of $e_1$, $e_2$, and
there is a $C^\infty$-function $K_{\phi,D}$
on $M^2\setminus \Sigma_\phi$ such that
\begin{equation}\label{eq:K}
    d\omega=K_{\phi,D}\, d\hat A.
\end{equation}
We call $K_{\phi,D}$ the {\it Gaussian curvature\/} of $D$ with respect
to $\phi$.
Let $\bar D$ be the pull-back of $D$ on $M^2\setminus \Sigma_{\phi}$.
Let $\sigma(t)$ be a regular curve on $U\setminus \Sigma_\phi$
with the arclength parameter $t$ with respect to
$ds^2=\phi^*\inner{~}{~}$.
We take a unit normal vector  $n(t)$ such that $(\dot \sigma,n)$
gives a positive frame on $TM^2$.
On the other hand, we take  $\hat n(t)\in E$
such that $(\phi(\dot \sigma),\hat n)$
gives a positive frame on $E$.
We can define two geodesic curvatures;
\[
  \kappa_g=ds^2(\bar D_t\dot\sigma(t),n(t)),
   \qquad
   \hat \kappa_g=\inner{D_t\phi(\dot\sigma(t))}{\hat n(t)}.
\]
Here, $\hat \kappa_g(t)$ is well-defined 
even when $\sigma(t)$ passes through the set $\Sigma_\phi$.
Since $\phi(n)=\sgn(\lambda) \hat n$, it holds that
$\kappa_g=\sgn(\lambda)\hat \kappa_g$.
We set $(\bar e_1,\bar e_2)=(\phi^{-1}(e_1),\phi^{-1}(e_2))$
if $U\subset M^+$ and set  
$(\bar e_1,\bar e_2)=(\phi^{-1}(e_2),\phi^{-1}(e_1))$
if $U\subset M^-$. 
Then $(\bar e_1,\bar e_2)$ gives an oriented orthonormal frame
on $TM^2$, 
and there is a $C^\infty$-function $\theta=\theta(t)$
such that
$\dot \sigma=\cos \theta \bar e_1+\sin \theta \bar e_2$
and 
$n=-\sin \theta \,\bar e_1+\cos \theta \,\bar e_2$.
Then we  get \begin{equation}\label{eq:stokes}
   \kappa_g dt=d\theta-(\sgn{\lambda})\omega.
\end{equation}
If the connection $D$ satisfies
the condition
\begin{equation} \label{eq:C}
   D_X\phi(Y)-D_Y\phi(X)-\phi([X,Y])=0
\end{equation}
for all vector fields $X,Y$ on $M^2$,
$(E,\inner{~}{~},D,\phi)$ is called a {\it coherent tangent bundle}.
Under the condition \eqref{eq:C},
$\bar D$ gives the Levi-Civita connection of $ds^2$ on 
$M^2\setminus \Sigma_\phi$, and $K_{\phi,D}$ coincides with 
the usual Gaussian curvature.
We consider a contractible 
triangular domain $\triangle \pt{ABC}$ on 
$M^2\setminus \Sigma_{\phi}$
such that it lies on the
left-hand side of the regular arcs $\pt{AB}$, $\pt{BC}$, $\pt{CA}$
which meet transversally at $\pt{A}$, $\pt{B}$, $\pt{C}\in M^2$.
By applying the Stokes formula, \eqref{eq:K} and \eqref{eq:stokes}
yield
that
\begin{equation}\label{eq:local}
    \angle\pt{A}+\angle\pt{B}
       +\angle\pt{C}-\pi=
    \int_{\partial \triangle\pt{ABC}} 
     \kappa_g \,d\tau + 
     \int_{\triangle\pt{ABC}} K_{\phi,D}\, dA,
\end{equation}
where $\angle\pt{A}$, $\angle\pt{B}$, $\angle\pt{C}$
are the interior angles of the domain $\triangle \pt{ABC}$.
To prove this, we do not need to assume that $\bar D$ is the Levi-Civita
connection.
However, we must remember that $K_{\phi,D}$ is not the usual Gaussian
curvature. 
One crucial point in this setting is that
\[
   \int_{M^2}K_{\phi,D}\, d\hat A=\frac1{2\pi}\int_{M^2}d\omega
\]
coincides with the Euler characteristic $\chi_{E}^{}$ 
of the vector bundle $E$.
In \cite{SUY2} (see also \cite{SUY1}), 
the authors gave the following two Gauss-Bonnet type formulas
\begin{equation}\label{eq:main}
  \chi_{E}^{}=
        \chi(M^+)-\chi(M^-)+S_+-S_-, 
  \qquad
  2\pi\chi(M^2)=
  \int_{M^2}K_{\phi,D}\, dA+
     2\int_{\Sigma_{\phi}} \kappa_s\, d\tau, 
\end{equation}
under the assumption that $(E,\inner{~}{~},D,\phi)$ is 
a coherent tangent bundle,
where $d\tau$ is the arclength element on the singular set
and $S_+,S_-$ are the numbers of positive and negative 
$A_3$-points, respectively.
After the publication of \cite{SUY2},
the authors found that the proof in \cite{SUY2} is
based only on the formula \eqref{eq:local} and the identity 
$\kappa_g=\sgn(\lambda)\hat \kappa_g$.
So we can conclude that the two formulas \eqref{eq:main}
hold without assuming \eqref{eq:C}.
Moreover, we can generalize these two formulas
to $\phi$ admitting more general singularities;
in other words, Theorem B in \cite{SUY2} holds
on $\phi$ without assuming \eqref{eq:C}.  
If $E=TM^2$, then $\chi_{E}^{}$
coincides with $\chi(M^2)=\chi(M^+)+\chi(M^-)$ in our setting.
So we get the following 
\begin{theorem}\label{theorem:gBW}
  Let $\phi:TM^2\to TM^2$ be a bundle homomorphism
  whose singular set 
  consists only of $A_2$ and $A_3$-points.
  Then 
  $2\chi(M^-)= S_+-S_-$
and $\int_{M^-}K_{\phi,D}\, d\hat A=\int_{\Sigma_\phi}\kappa_s\,d\tau$
hold.
\end{theorem}
Let  $f:M^2\to (N^3,g)$ be an  immersion into
an orientable Riemannian $3$-manifold.
Then there is a globally defined unit normal
vector field $\nu$ along $f$.
We define the shape operator
$\phi:TM^2\ni v \mapsto -D_v\nu\in TM^2$,
as a bundle homomorphism, where $D$
is the Levi-Civita connection of $(N^3,g)$.
A singular point of $\phi$ is called 
an {\it inflection point\/} of $f$.
We get the following
\begin{corollary}[A generalization of the Bleeker-Wilson formula]%
\label{cor:gBW0}
 Suppose that the shape operator admits only $A_2$ and
 $A_3$-points.
 Then 
 $2\chi(M^-)=I_+-I_-$ holds, 
 where $I_+$ {\rm(}resp.\ $I_-${\rm)} is 
 the number of positive {\rm(}resp.\ negative{\rm)}
 $A_3$-inflection points.
\end{corollary}
The original formula was for the case  $N^3=\R^3$
(see \cite{BW}).
In \cite{SUY5}, the authors pointed out that  the formula holds for
space forms.
Also, they gave in \cite{SUY5} several applications
of \eqref{eq:main} under the assumption \eqref{eq:C}.  
However, now we can remove \eqref{eq:C},
and we get also the results that follow here. 

\section{Rotation of vector fields.}
We fix a Riemannian metric $ds^2$ on $M^2$.
There is a unique 2-form $\mu$ on $M^2$
such that $\mu(e_1,e_2)=1$ where $e_1$, $e_2$
is a local oriented orthonormal frame field
on $M^2$.
Let $X$ be a vector field on $M^2$.
The $C^\infty$-function
$\rot(X):=\mu(D_{e_1}X,D_{e_2}X)$
defined on $M^2$
is called {\it the rotation\/} of $X$,
where $D$ is the Levi-Civita connection of $(M^2,ds^2)$.
Consider a bundle homomorphism
$\phi:TM^2\ni v \mapsto D_vX\in TM^2$.
The singular set $\Sigma_X$ of $\phi$ coincides
with the zeros of $\rot(X)$, 
called the set of {\it irrotational points}.
Moreover, an $A_3$-singular point
is called an {\it irrotational cusp}.
In fact, if $M^2=\R^2$ is the Euclidean plane,
then $X$ induces a
map $\tilde X:\R^2\to \R^2$, and $A_3$ (resp. $A_2$)
points correspond to cusps (resp.\ folds)
of $\tilde X$ (see \cite{SUY5}).
Suppose that
$X$ admits only $A_2$ and
$A_3$-irrotational points.
Then $\Sigma_X$
consists of a finite disjoint union of closed
regular curves 
$\gamma_1$, \dots,$\gamma_m$
on $M^2$ such that $M^+$ lies in the left hand side of
each $\gamma_j$.
Then the singular curvature on $\gamma_j$
is given by
$\kappa_s:=\mu(\dot X,\ddot X)/|\dot X|^3$
(we propose to call it the {\it irrotational curvature}),
where $\dot X=D_{\dot \gamma_j(t)}X$ and 
$\ddot X=D_{\dot \gamma_j(t)}\dot X$.
The following assertion follows directly from Theorem \ref{theorem:gBW}.
\begin{e-proposition} 
 Suppose that $X$ admits only $A_2$ and
 $A_3$-irrotational points.
 Then 
 it holds that
 \begin{multline*}
     2\chi(M^-)=
        C_+-C_-, \quad
        \int_{M^-}K_{\phi,D} d\hat A=
        \int_{\Sigma_X}\kappa_s d\tau,
      \\
      M^-:=\biggl\{p\in M^2\,;\,
       \rot(X)_p<0\biggr\},
 \end{multline*}
 where $C_+$ {\rm(}resp.\ $C_-${\rm)} is 
 the number of positive {\rm(}resp.\ negative{\rm)}
 irrotational cusps.
\end{e-proposition}

\section{Singularities of Blaschke normal maps on convex surfaces.}
Let $S^2$ be a $2$-sphere and
$f:S^2\to \R^3$ a strictly convex embedding.
In affine differential geometry,
it is well-known that there are a transversal vector field $\xi$ along
$f$, a torsion free connection $\nabla$, a bundle homomorphism
$\alpha:TS^2\to TS^2$ (called the {\it affine shape operator}), 
and a positive definite symmetric covariant tensor $h$ such that 
(cf.\ \cite{NS})
$D_XY=\nabla_XY+h(X,Y)\xi$ and
$D_X\xi=-\alpha(X)$ for any vector fields $X$, $Y$ 
on $S^2$, where
$D$ is the canonical affine connection on $\R^3$.
Moreover, such a structure $(\xi,\nabla,\alpha,h)$
is uniquely determined up to a constant multiplication of $\xi$.
Here $\xi$ induces a map $\tilde \xi:S^2\to \R^3$ called 
the {\it Blaschke normal map}.
It is obvious that the singular points of $\alpha$ coincides 
with those of $\tilde \xi$.
\begin{lemma}
 The  Blaschke normal map $\tilde \xi$ is
 a wave front $($cf. {\rm \cite{AGV}} for 
the definition of wave front$)$.
\end{lemma}
\begin{proof}
 Consider a non-zero section 
$L:S^2\ni p \mapsto (\tilde \xi_p,\nu_p)
\in T^*\R^3=\R^3\times (\R^3)^*$,
 where $\nu:S^2\to (\R^3)^*$ is the map into 
 the dual vector space $(\R^3)^*$ of $\R^3$
 such that $\nu_p(\tilde\xi_p)=1$ 
 and $\nu_p(df(T_pS^2))=\{0\}$ for each $p\in S^2$.
 Take a local coordinate system $(u_1,u_2)$
 of $S^2$.
 Then we have that
 \begin{multline*}
    \nu_{u_i}(f_{u_j})=
      D_{\partial_i}\nu(f_{u_j})=\nu(D_{\partial_i}f_{u_j})\\
      =-\nu\left(\nabla_{\partial_i} \partial_j
                +h(\partial_i,\partial_j)\tilde \xi\right)
        =-h(\partial_i,\partial_j) \quad (i,j=1,2),
 \end{multline*}
 where $\partial_i:=\partial /\partial u_i$
 and $f_{u_i}:=df(\partial_i)$.
 Since $h$ is positive definite, 
 $\nu_{u_1}$, $\nu_{u_2}$ are linearly independent. 
 Moreover, $\nu,\nu_{u_1},\nu_{u_2}$ are also 
 linearly independent, since $\nu(df(T_pS^2))=0$.
 In particular, $L$ induces
a Legendrian immersion of
$S^2$ into the projective cotangent 
 bundle $P(T^*\R^3)$
of $T^*\R^3$. 
\end{proof}

By applying the criteria of cuspidal edges and 
swallowtails (cf. \cite{SUY5}, $A_2$ and $A_3$-points
correspond to the cuspidal edges and
swallowtails of the Blaschke normal map $\tilde \xi$.
So we get the following
\begin{theorem} 
 Suppose that  $\tilde \xi$ admits only cuspidal edges and swallowtails.
 Then 
 $2\chi(M^-)=S_+-S_-$ holds,
 where $M^-:=\{p\in S^2\,;\, \det(\alpha(p))<0\}$ and $S_+$ {\rm(}resp. $S_-${\rm)} is 
 the number of positive {\rm(}resp. negative{\rm)}
 swallowtails of $\tilde \xi$.
\end{theorem}

A different formula for $S_++S_-$ is given by Izumiya-Marer
\cite{IM}.

\end{document}